\documentclass{elsart3-1}

\usepackage{graphicx}
\usepackage{amssymb}

\usepackage[english,francais]{babel}

\newtheorem{theorem}{Theorem}[section]

\newtheorem{e-proposition}[theorem]{Proposition}

\newtheorem{e-definition}[theorem]{Definition\rm}

\renewcommand{\theequation}{\arabic{equation}}
\def\og{\leavevmode\raise.3ex\hbox{$\scriptscriptstyle\langle\!\langle$~}}
\def\fg{\leavevmode\raise.3ex\hbox{~$\!\scriptscriptstyle\,\rangle\!\rangle$}}

\journal{the Acad\'emie des sciences}
\begin{document}
% place in the next line the header (rubrique) chosen for your article,
% if you know it (you can also have 2, format : Header1/Header2
\centerline{}
\begin{frontmatter}

% Title, authors and addresses

% use the thanksref command within \title, \author or \address for footnotes;
% use the ead command for the email address,
% and the form \ead[url] for the home page:
% \title{Title\thanksref{label1}}
% \thanks[label1]{}
% \author{Name\thanksref{label2}}
% \ead{email address}
% \ead[url]{home page}
% \thanks[label2]{}
% \address{Address\thanksref{label3}}
% \thanks[label3]{}
\selectlanguage{english}
\title{ New Extended Formulations of Euler-Korteweg Equations \\
          Based on a  Generalization of the Quantum Bohm Identity}

% use optional labels to link authors explicitly to addresses:
% \author[label1,label2]{}
% \address[label1]{}
% \address[label2]{}
% The [label1] can be suppressed if there is only one address for all authors

\selectlanguage{english}
\author[authorlabel1]{Didier Bresch\thanksref{lab1}},
\ead{Didier.Bresch@univ-savoie.fr}
\author{Fr\'ed\'eric Couderc},
\ead{frederic.couderc@math.univ-toulouse.fr}
\author{Pascal Noble\thanksref{lab2}},
\ead{pascal.noble@math.univ-toulouse.fr}
\author[authorlabel2]{Jean--Paul Vila}
\ead{vila@insa-toulouse.fr}
\thanks[lab1]{Research of D.B. was partially supported by the ANR project DYFICOLTI ANR-13-BS01-0003- 01}
\thanks[lab2]{Research of P.N. was partially supported by the ANR project BoND ANR-13-BS01-0009-01 .}
\address[authorlabel1]{LAMA -- UMR5127 CNRS, Bat. Le Chablais, Campus Scientifique, 73376 Le Bourget du Lac, France}
\address[authorlabel2]{IMT, INSA Toulouse, 135 avenue de Rangueil, 31077 Toulouse Cedex 9, France}
%\address[authorlabel3]{Address3}

% If you know the dates of reception, and acceptation you can put them now;
%  idem the name of the person presenting the Note

\medskip
\begin{center}
{\small Received *****; accepted after revision +++++\\
Presented by £££££}
\end{center}

\begin{abstract}
\selectlanguage{english}
% Text of abstract in English
In this note, we  propose an original  extended formulation of Euler-Korteweg systems based
on  a generalization of the quantum Bohm potential identity.  This new formulation allows
to propose a useful  construction of a numerical scheme  with entropy stability property  
under a hyperbolic CFL condition.  We also comment  the use of the identity for
compressible Navier-Stokes equations with degenerate viscosities.
%PN{\it To cite this article: A. Name1, A. Name2, C. R. Acad. Sci. Paris, Ser. I 340 (2005).}

\selectlanguage{francais}
% Text of abstract in French
\noindent{\bf R\'esum\'e} \vskip 0.5\baselineskip \noindent
{\bf G\'en\'eralisation de l'identit\'e de Bohm quantique et nouvelles formulations 
augment\'ees pour  \'equations d'Euler-Korteweg. }
 Dans cette note, on propose une formulation augment\'ee originale
 des syst\`emes d'Euler-Korteweg bas\'ee sur une g\'en\'eralisation de 
 l'identit\'e  dite du potentiel de Bohm quantique . La motivation principale est la 
 construction de sch\'emas avec stabilit\'e entropique sous condition CFL
 hyperbolique du syst\`eme d'Euler-Korteweg.  On pr\'esente \'egalement 
 quelques commentaires concernant l'utilisation de cette identit\'e dans le
 cadre des \'equations de Navier--Stokes avec 
 viscosit\'es d\'eg\'en\'er\'ees.

%PN{\it Pour citer cet article~: A. Name1, A. Name2, C. R. Acad. Sci. Paris, Ser. I 340 (2005).}
\end{abstract}
\end{frontmatter}

% now the Version française abrégée, if it exists
\selectlanguage{francais}
\section*{Version fran\c{c}aise abr\'eg\'ee}
Dans cette note, nous introduisons de nouvelles formulations augment\'ees (sous forme conservative) du syst\`eme  d'Euler-Korteweg (\ref{EKma})--(\ref{EKCa}) en plusieurs dimensions d'espace. Selon le choix des coefficients de capillarit\'e, ce type de syst\`eme intervient dans la mod\'elisation 
des m\'elanges de type liquide-vapeur, des super-fluides ou de l'hydrodynamique quantique 
par exemple. 
   Nous nous int\'eressons ici \`a de nouvelles formulations permettant de construire un
sch\'ema num\'erique \`a stabilit\'e entropique. Nous \'etendons ainsi les travaux r\'ecents
de deux des auteurs qui consid\'eraient le cas uni-dimensionnel. Nous traitons
les termes dispersifs implicitement et donnons un r\'esultat de stabilit\'e entropique des sch\'emas
d'ordre 1 sous condition CFL hyperbolique. A titre d'illustration, nous pr\'esentons
des r\'esultats num\'eriques pour des films minces avec tension de surface mod\'elis\'es
par les \'equations de Saint-Venant.  

\selectlanguage{english}
% main text
\section{Introduction}
\label{sec1}
In this paper, we introduce (new) extended formulations of the so-called Euler Korteweg system which, in several space dimensions, reads
{\setlength\arraycolsep{1pt}
\begin{eqnarray}
\label{EKma}
\displaystyle
\partial_t\varrho&+&{\rm div}(\varrho{\bf u})=0,\\
\label{EKmo}
\displaystyle
\partial_t (\varrho{\bf u})&+& {\rm div}(\varrho  {\bf u}\otimes{\bf u}) +\nabla p(\varrho)={\rm div}(\mathbf{K}),
\end{eqnarray}
}
\noindent
where $\varrho$ denotes the fluid density, ${\bf u}$ the fluid velocity, $p(\varrho)$ the fluid pressure and $\mathbf{K}$ the Korteweg stress tensor defined as
\begin{equation} \label{EKCa}
\displaystyle
\mathbf{K}=\left(\varrho{\rm div}(K(\varrho)\nabla \varrho)+\frac{1}{2}(K(\varrho)-\varrho K'(\varrho))|\nabla \varrho|^2 \right)\mathbf{I}_{\mathbb{R}^n}-K(\varrho)\nabla\varrho\otimes \nabla \varrho.
\end{equation}
\noindent
with $K(\varrho)$ the capillary coefficient. These models comprise liquid-vapor mixture (for instance highly pressurized and hot water in nuclear reactors cooling system) \cite{JTB}, superfluids (Helium near absolute zero) \cite{HAC} or even regular fluids at sufficiently small scales (think of ripples on shallow waters) \cite{LiScGo}. In quantum hydrodynamic, the capillary coefficient is chosen so that $\varrho K(\varrho)=constant$: in this case; the Euler-Korteweg equations correspond to the nonlinear Schr\"odinger equation after Madelung transform. In classical fluid mechanics, the capillary coefficient $K(\varrho)$ is chosen constant.\\

\noindent
The system (\ref{EKma})--(\ref{EKmo}) admits two additional conservations laws. One conservation law is satisfied by the fluid velocity 
\begin{equation}\label{EKfl}
\displaystyle
\partial_t {\bf u}+{\bf u}\cdot \nabla\,{\bf u}+\nabla(\delta\mathcal{E})=0,
\end{equation}
\noindent
with $\mathcal{E}$ the potential energy and $\delta\mathcal{E}$ its {\it variational} gradient
\begin{equation}\label{defEN}
\displaystyle
\mathcal{E}(\varrho,\nabla\varrho)=F_0(\varrho)+\frac{1}{2}K(\varrho)|\nabla\varrho|^2,\quad \delta\mathcal{E}=F'_0(\varrho)-\frac{1}{2}K'(\varrho)|\nabla\varrho|^2-K(\varrho)\Delta\varrho.
\end{equation}
\noindent
The local existence of strong solution to (\ref{EKma}), (\ref{EKfl}) is proved in \cite{BeDaDe}. For that purpose, the authors introduced an extended formulation by considering an additional velocity ${\bf w}=\nabla\varphi(\varrho)$ with $\sqrt{\varrho}\varphi'(\varrho)=\sqrt{K(\varrho)}$: 
%{\setlength\arraycolsep{1pt}
\begin{equation}
\displaystyle
\partial_t{\bf u}+{\bf u}\cdot \nabla\,{\bf u}+\nabla\left(F'_0(\varrho)-\frac{|{\bf w}|^2}{2}\right)=\nabla\left(a(\varrho){\rm div}({\bf w})\right),\quad
%\displaystyle
\partial_t{\bf w}+\nabla\left({\bf u}^T\,{\bf w}\right)=-\nabla\left(a(\varrho){\rm div}({\bf u})\right),\nonumber
\end{equation}
%}
\noindent
where $a(\varrho)=\sqrt{K(\varrho)\varrho}$. This formulation is particularly adapted to the derivation of a priori estimates, the very first one being a conservation law on the total (kinematic+potential) energy :
{\setlength\arraycolsep{1pt}
\begin{eqnarray}
\displaystyle
\partial_t\left(\frac{\varrho}{2}|{\bf u}|^2+\mathcal{E}(\varrho,\nabla\varrho)\right)
+{\rm div}\left({\bf u}\left(\frac{\varrho}{2}|{\bf u}|^2+\mathcal{E}(\varrho,\nabla\varrho)+p(\varrho)\right)\right)=
{\rm div}\Bigl(F(\varrho)(\nabla{\bf w}\,{\bf u}-\nabla{\bf u}\,{\bf w})\Bigr)&&\nonumber\\
\label{EKen}
-{\rm div} \Bigl((F(\varrho)-\varrho F'(\varrho))({\rm div}({\bf w}){\bf u}-{\rm div}({\bf u}){\bf w})\Bigr).&&
\end{eqnarray}
}
with $F'(\varrho)=\varrho\varphi'(\varrho)$.\\

\noindent
In this note, we introduce a new extended formulation of (\ref{EKma})--(\ref{EKmo}) by considering the conservative variables $\varrho{\bf u}, \varrho{\bf w}$ instead of ${\bf u}, {\bf w}$.  The key point is
the generalization of the quantum potential Bohm  identity. It allows to transform the Euler-Korteweg system into a hyperbolic system perturbed by a second order skew symmetric term. The main motivation is the construction of a numerical scheme which is easily proved ``entropy'' stable. 
Two of the authors performed this approach in a one dimensional setting and proved entropy stability under ``capillary'' Courant-Friedrichs-Lewy (denoted CFL in the sequel) condition  \cite{NoVi}. Here, we extend this approach to the multi-dimensional setting. Moreover, in order to avoid restrictive CFL condition, we treat dispersive terms implicitely and prove entropy stability of first order schemes under a hyperbolic CFL condition. We present preliminary numerical results for thin films with surface tension modeled by the shallow water equations.

\section{Generalization of the quantum Bohm identity and extended formulation}
\label{sec2}

Let us first present an extension of the quantum potential Bohm identity
$$2 \varrho {\rm div}\bigl(\Delta \sqrt{\varrho}/\sqrt{\varrho}\bigr)
   = {\rm div}\Bigl(\varrho\nabla\nabla \log \varrho\Bigr)$$
 strongly used in quantum fluid mechanics.
More precisely, we can prove after some algebraic calculations the following relation
\begin{equation} \label{Bohm}
\varrho \nabla \bigl(\sqrt{K(\varrho)}\, \Delta( \int_0^\varrho \sqrt{K(s)}\, ds)\bigr)
   = {\rm div}\Bigl(F(\varrho) \nabla\nabla \varphi(\varrho)\Bigr)
   - \nabla \Bigl(\bigl(F(\varrho)-F'(\varrho)\varrho\bigr) \Delta \varphi(\varrho)\Bigr)
\end{equation}
with $\sqrt{\varrho} \varphi'(\varrho) = \sqrt{K(\varrho)}$, $F'(\varrho)= \sqrt{K(\varrho)\,\varrho}$.
This relation is a non trivial extension of the quantum Bohm identity which corresponds to the case  $K(\varrho)=c/\varrho$. Remark that  the left--hand side of  (\ref{Bohm}) corresponds to the capillarity
term $ K(\varrho) |\nabla\varrho|^2/2$. It suffices to observe that
$$ K(\varrho) |\nabla\varrho|^2  =  |\nabla( \int_0^\varrho \sqrt{K(s)})\, ds|^2 $$
and  thus as observed in \cite{BrDeLi} the variational gradient of the potential energy may be written
$$\delta \mathcal{E}= F'_0(\varrho)
     - \bigl(\sqrt{K(\varrho)}\, \Delta( \int_0^\varrho \sqrt{K(s)}\, ds)\bigr).$$

\noindent
Following now the strategy of \cite{BeDaDe}, we introduce a ``good'' additional unknown, homogeneous to a velocity.  We denote this additional velocity ${\bf w}=\nabla\varphi(\varrho)$ with $\sqrt{\varrho}\varphi'(\varrho)=\displaystyle\sqrt{K(\varrho)}$.
In order to write a suitable extended formulation of the Euler Korteweg model, we also define $F(\varrho)$ so that
$F'(\varrho)=\sqrt{K(\varrho)\varrho}$. The Euler Korteweg system admits the extended formulation
\begin{equation}\label{EKm}
\left\{
\begin{array}{l}
\vspace{0.2cm}
\partial_t\varrho+{\rm div}(\varrho\,{\bf u})=0,\\
\vspace{0.2cm}
\partial_t(\varrho\,{\bf u})+{\rm div}(\varrho  {\bf u}\otimes{\bf u}) +\nabla p(\varrho)={\rm div}(F(\varrho)\nabla {\bf w}^T)-\nabla\left((F(\varrho)-\varrho F'(\varrho)){\rm div}({\bf w})\right),\\
\vspace{0.2cm}
\partial_t(\varrho\,{\bf w})+ {\rm div}(\varrho  {\bf w}\otimes{\bf u}) =-{\rm div}(F(\varrho)\nabla {\bf u}^T)+\nabla\left((F(\varrho)-\varrho F'(\varrho)){\rm div}({\bf u})\right),
\end{array}\right.
\end{equation}

\noindent
The total energy is then transformed into a classical entropy of the first order part of (\ref{EKm}) 
$$
\displaystyle
\frac{\varrho}{2}\|{\bf u}\|^2+\mathcal{E}(\varrho, \nabla\varrho)=\frac{\varrho}{2}\left(\|{\bf u}\|^2+\|{\bf w}\|^2\right)+F_0(\varrho):={\bar\mathcal{E}}(\varrho, {\bf u}, {\bf w}),
$$
\noindent
whereas the second order part is skew symmetric. As a consequence, the energy conservation law (\ref{EKen}) is obtained through a similar computation than in the first order case.

\bigskip

\noindent {\bf Remark.}
Note that our formulation may be coupled to the result recently obtained in \cite{BrDeZa} to
write an augmented formulation to the following compressible Navier--Stokes system
with drag and capillary terms
{\setlength\arraycolsep{1pt}
\begin{eqnarray}
\label{EKma1}
\displaystyle
\partial_t\varrho&+&{\rm div}(\varrho{\bf u})=0,\\
\label{EKmo1}
\displaystyle
\partial_t (\varrho{\bf u})&+& {\rm div}(\varrho  {\bf u}\otimes{\bf u}) 
     - 2 {\rm div}(\mu(\varrho) D({\bf u}) - \nabla(\lambda(\varrho) {\rm div} {\bf u}) 
      + \nabla p(\varrho) + r_1 \varrho |u|^\alpha u ={\rm div}(\mathbf{K}),
\end{eqnarray}
}
if  $\lambda(\varrho) = 2(\mu'(\varrho)\varrho - \mu(\varrho))$ and 
$K(\varrho)= c (\mu'(\varrho))^2/\varrho$.  
Such compatible system could be used to prove global existence of weak solutions to the compressible Navier-Stokes equations with degenerate viscosities without capillary and drag terms.
  It has been recently performed in \cite{VaYu} in the case $\mu(\varrho)= \mu \varrho$ and  $\lambda(\varrho)=0$ introducing the quantum capillary term namely with  $K(\varrho)=c/\varrho$.
   
\section{Application: stable schemes under hyperbolic CFL condition}
\label{sec3}

In this section, we introduce a numerical scheme for (\ref{EKm}). The numerical domain is a rectangle defined by $0\leq x\leq L_x$ and $0\leq y\leq L_y$, which is divided into $N=n_x\times n_y$ rectangular cells. For the sake of simplicity, we consider uniform grid with constant spatial steps $\delta x$ and $\delta y$. We focus on the spatial discretization of the second order terms: they are written as
$\displaystyle \partial_{\alpha}(f(\varrho)\partial_\beta u)$ with $(\alpha,\beta)\in\{x, y\}$ and $f(\varrho)=F(\varrho), \varrho F'(\varrho),F(\varrho)-\varrho F'(\varrho)$. For that purpose, we introduce the following finite difference operators:
\begin{equation}\label{discrete-diff-op}
\begin{array}{ll}
\displaystyle
(d_1 u)_{i,j}=\frac{u_{i+1/2,j}-u_{i-1/2,j}}{\delta x},\quad (d_1^+u)_{i+1/2,j}=\frac{u_{i+1,j}-u_{i,j}}{\delta x},\quad ({\bar d}_1u)_{i,j}=\frac{u_{i+1,j}-u_{i-1,j}}{2\delta x},\\
\displaystyle
(d_2 u)_{i,j}=\frac{u_{i,j+1/2}-u_{i,j-1/2}}{\delta y},\quad (d_2^+u)_{i,j+1/2}=\frac{u_{i,j+1}-u_{i,j}}{\delta y},\quad ({\bar d}_2 u)_{i,j}=\frac{u_{i,j+1}-u_{i,j-1}}{2\delta y}.
\end{array}
\end{equation}
\noindent
As a result the differential operator $\displaystyle \mathcal{T}(\varrho){\bf u}={\rm div}(F(\varrho)\nabla{\bf u}^T)+\nabla\left((\varrho F'(\varrho)-F(\varrho)){\rm div}({\bf u})\right)$ is approximated by $\mathcal{T}_h(\varrho)$ defined as
$$
\mathcal{T}_h(\varrho){\bf u}_{i,j}=\left\{\begin{array}{ll}
\displaystyle
d_1\left(\varrho\,F'(\varrho)d_1^+{\bf u}_1\right)_{i,j}+{\bar d}_2\left(F(\varrho){\bar d}_1{\bf u}_2\right)_{i,j}+{\bar d}_1\left((\varrho F'(\varrho)-F(\varrho)){\bar d}_2{\bf u}_2\right)_{i,j},\\
\displaystyle
{\bar d}_1\left(F(\varrho){\bar d}_2{\bf u}_1\right)_{i,j}+{\bar d}_2\left((\varrho F'(\varrho)-F(\varrho)){\bar d}_1{\bf u}_1\right)_{i,j}+d_2\left(\varrho F'(\varrho)d_2^+{\bf u}_2\right)_{i,j}
\end{array}\right.
$$
\noindent
We discretize (\ref{EKm}) as follows
\begin{equation}\label{EKd}
\left\{
\begin{array}{l}
\vspace{0.2cm}
\frac{\varrho_{i,j}^{n+1}-\varrho_{i,j}^n}{\delta t}+d_1(\mathcal{F}_{\varrho,1}^n)_{i,j}+d_2(\mathcal{F}_{\varrho,2}^n)_{i,j}=0,\\
\vspace{0.2cm}
\frac{(\varrho\,{\bf u})_{i,j}^{n+1}-(\varrho{\bf u})_{i,j}^n}{\delta t}+d_1({\bf \mathcal{F}}_{{\bf u},1}^n)_{i,j}+d_2({\bf \mathcal{F}}_{{\bf u},2}^n)_{i,j}=\mathcal{T}_h(\varrho^{n+1}){\bf w}^{n+1}_{i,j},\\
\vspace{0.2cm}
\frac{(\varrho{\bf w})_{i,j}^{n+1}-(\varrho{\bf w})_{i,j}^n}{\delta t}+d_1({\bf \mathcal{F}}_{{\bf w},1}^n)_{i,j}+d_2({\bf \mathcal{F}}_{{\bf w},2}^n)_{i,j} =-\mathcal{T}_h(\varrho^{n+1}){\bf u}^{n+1}_{i,j},
\end{array}\right.
\end{equation}
\noindent
where $\mathcal{F}_{\varrho,k}^n , {\bf \mathcal{F}}_{{\bf u},k}^n, {\bf \mathcal{F}}_{{\bf w},k}^n$ ($k=1,2$) are classical Rusanov fluxes evaluated  at $\varrho^n, {\bf u}^n, {\bf w}^n$. More precisely, the convection part is treated explicitly whereas the capillary terms are treated implicitly. Remark that there are no capillary terms in the mass conservation law so that the implicit step amounts to solve a {\it linear} sparse system and is easily proved entropy stable. As a consequence, one can prove, by using discrete duality properties of the discrete second order  operators, the following entropy stability result.\\

\begin{theorem}
Suppose (\ref{EKd}) is completed with periodic boundary conditions. Assume the hyperbolic scheme (system (\ref{EKd}) with $F=0$) is entropy stable then the fully hyperbolic/capillary scheme (\ref{EKd}) is entropy stable: 
$$
\displaystyle
\sum_{i=0}^{n_x}\sum_{j=0}^{n_y}{\bar\mathcal{E}}(\varrho_{i,j}^{n+1}, {\bf u}_{i,j}^{n+1}, {\bf w}_{i,j}^{n+1})\leq\sum_{i=0}^{n_x}\sum_{j=0}^{n_y}{\bar\mathcal{E}}(\varrho_{i,j}^{n}, {\bf u}_{i,j}^{n}, {\bf w}_{i,j}^{n}).
$$
\end{theorem}
\noindent
This means that the numerical scheme (\ref{EKd}) is entropy stable under a classical hyperbolic Courant-Friedrichs-Lewy condition. As an application, we carried out a numerical simulation of a thin film falling down an inclined plane. A consistent shallow water model \cite{BoChNoVi} is given by
{
\setlength\arraycolsep{1pt}
\begin{eqnarray}
\label{consmass}
\displaystyle
\partial_t h+{\rm div}\,(h\mathbf{u}\,)&=&0,\\
\displaystyle
\label{consmomentum}
\partial_t (h\mathbf{u})+{\rm div}\left(h\,\mathbf{u}\otimes\,\mathbf{u}\right)+\nabla(p(h))+\left(\frac{g\sin(\theta)}{\nu}\right)^2\,\partial_x\left(\frac{2h^5}{225}\right)\mathbf{e}_1&=&S(h, {\bf u})+\frac{\sigma h}{\rho}\nabla(\Delta h).
\end{eqnarray}
}

\noindent
with $p(h)=\displaystyle g\cos(\theta)h^2/2$ and $\displaystyle S(h, {\bf u})=gh\sin(\theta)\mathbf{e}_1-3{\nu\mathbf{u}}/{h}$ and ${\bf e}_1$ the first vector of the canonical base directed downstream. Here $g=9.8$ is the gravity constant, $\rho, \nu, \sigma$ are respectively the fluid density, kinematic viscosity and surface tension whereas $\theta$ is the inclination of the plane. We picked the values found in \cite{LiScGo} for a solution with $31\%$ glycerin by weight: $\rho=1.07\times 10^{3}\, kg.m^{-3}$, $\nu=2.3\times 10^{-6}\, m^2.s^{-1}$ and $\sigma=67\times 10^{-3}\, kg.s^{-2}$. The source term is treated implicitly: since the source term is only in the equation for ${\bf u}$ and is linear with respect to ${\bf u}$, the implicit step remains linear. We first carry out a numerical simulation of the original experience in \cite{LiScGo} but imposed periodic boundary conditions in both directions.\\

\begin{figure}[h!]
\begin{center}
\includegraphics[width=8cm]{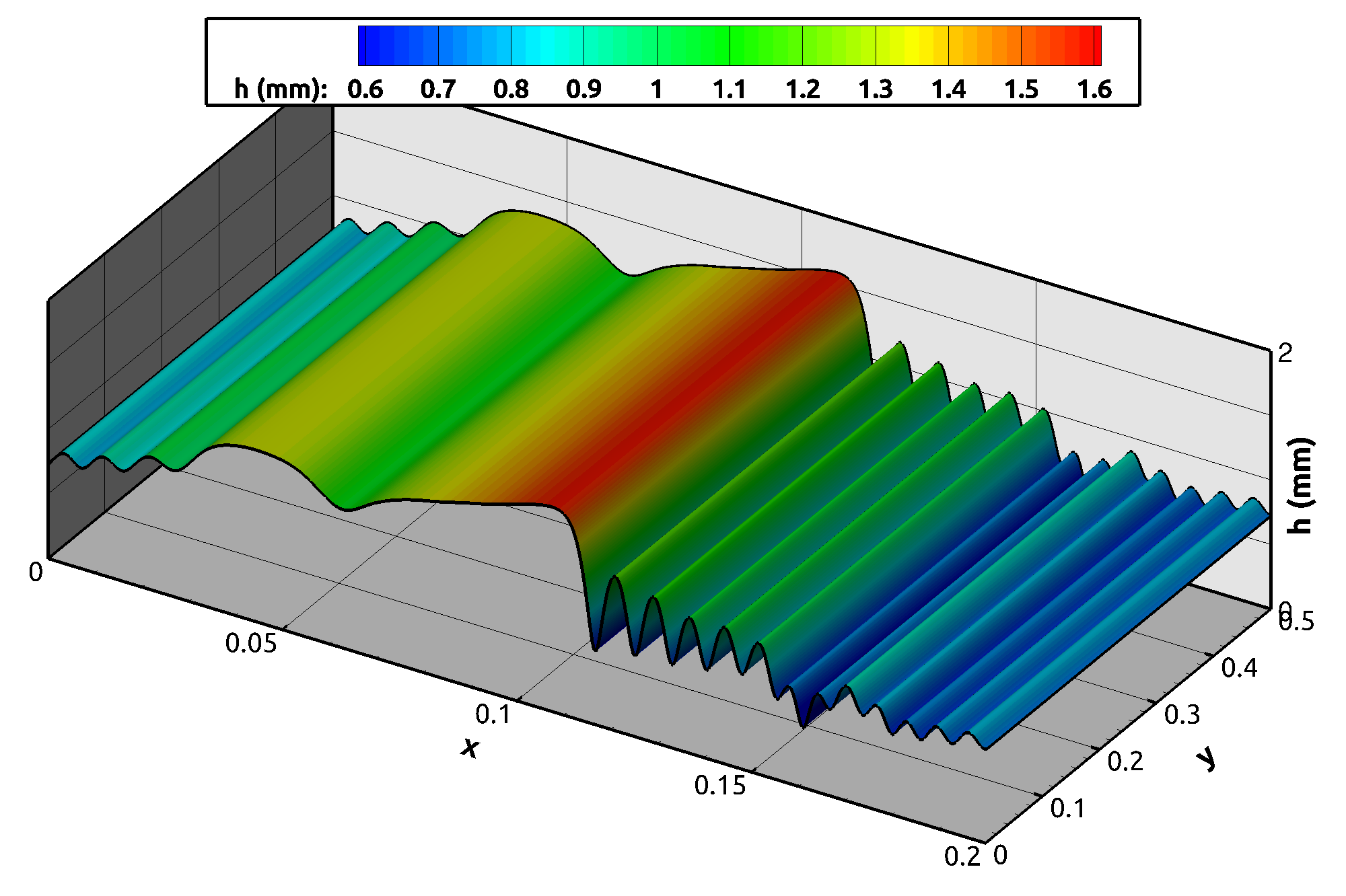}\hfill
\includegraphics[width=8cm]{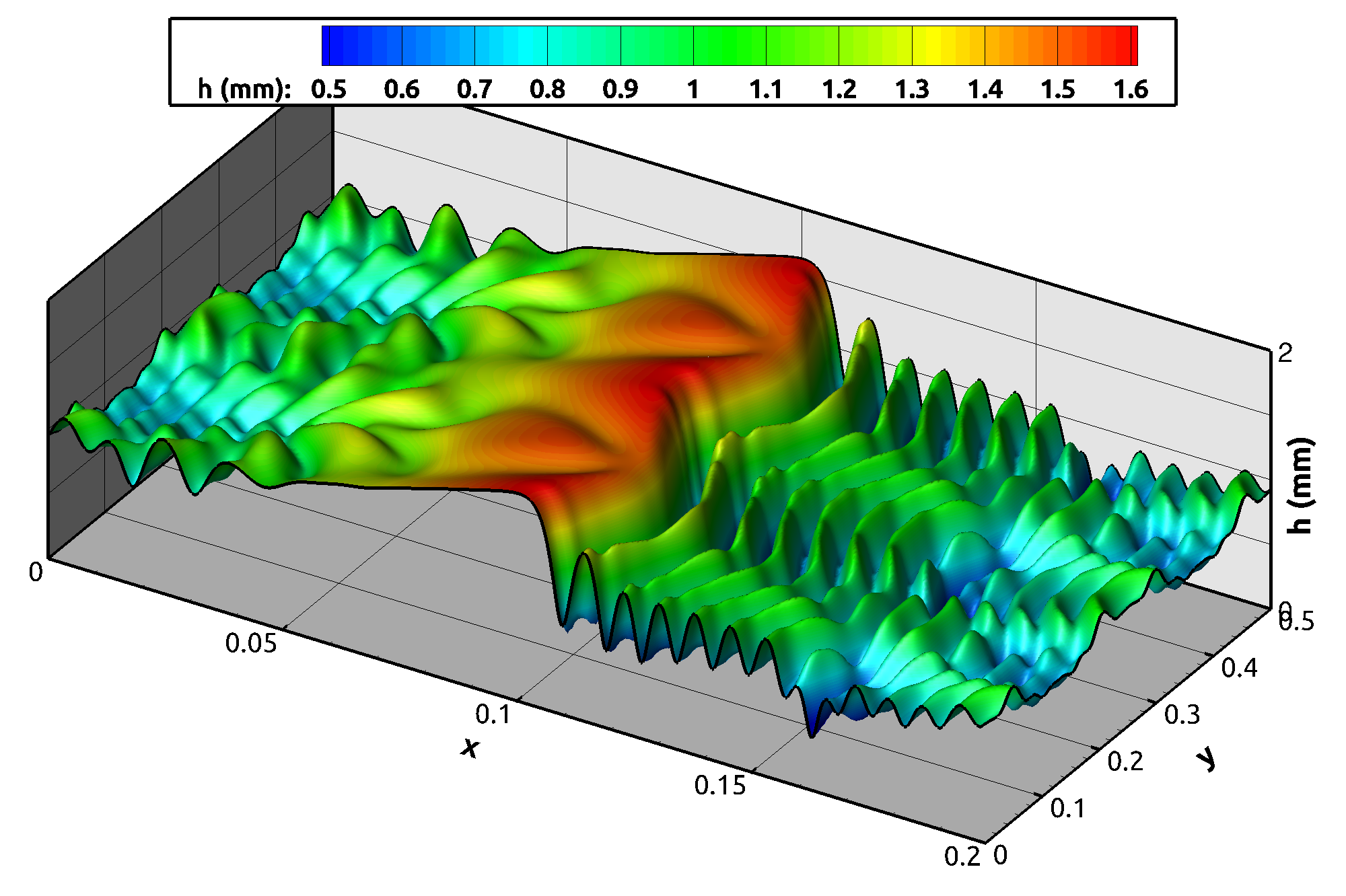}\hfill
\caption{Numerical simulation of a roll-wave in presence of surface tension. On the left: one dimensional roll-wave without transverse perturbations. On the right: a two-dimensional roll-wave}
\end{center}
\end{figure}

\noindent
In order to test the robusness of the scheme, we also carried various numerical experiments of a drop falling down a plane in order to deal with wet/dry fronts. For that purpose, we introduced a precusor film with a thickness  of $1.0\times 10^{-5}\,mm$.

\noindent
\begin{figure}[h!]
\begin{center}
\includegraphics[width=9cm]{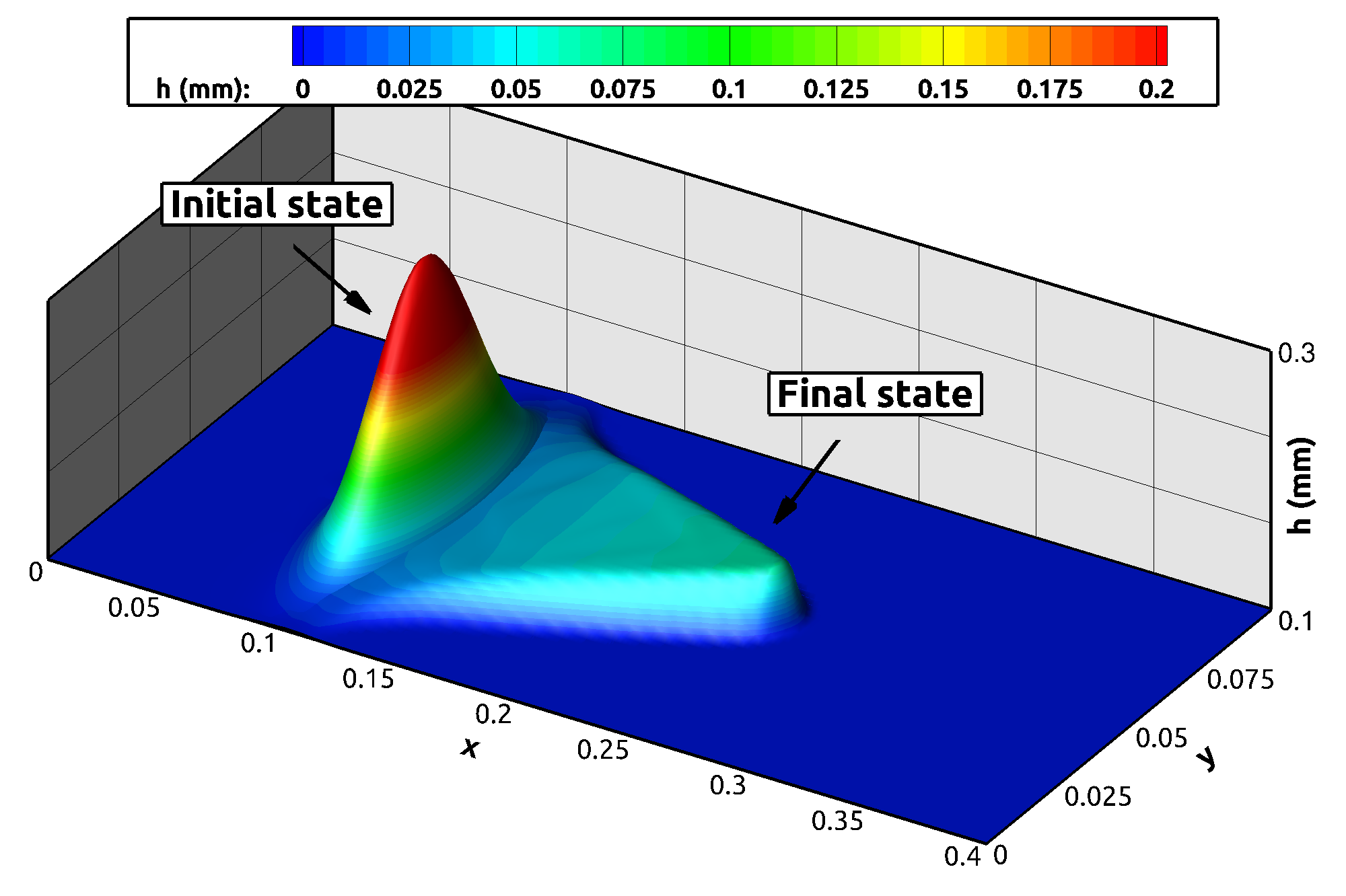}
\caption{Drop falling down an incline plane ($\theta=60^o$) at time $t=0$ and $t=1s$.  The fluid density, kinematic viscosity and surface tension are respectively $\rho=1.0\times 10^{3}\, kg.m^{-3}$, $\nu=1.0\times 10^{-6}\, m^2.s^{-1}$ and $\sigma=67\times 10^{-3}\, kg.s^{-2}$.}\label{fig2}
\end{center}
\end{figure}

\noindent
We will deal the problems of considering physical boundary conditions, deriving higher order schemes and considering wet/dry fronts in a forthcoming paper \cite{BrNoVi-forth}. This will be useful to compute instabilities in moving contact lines.

% etc, etc

% The Appendices part is started with the command \appendix;
% appendix sections are then done as normal sections
% \appendix

% \section{}
% \label{}

% The Acknowledgements are an un-numbered section
%\section*{Acknowledgements}
% Acknowledgements text here

\end{document}